\documentclass[12pt]{article}
\usepackage{amssymb}
\begin{document}

\author{David Carf\`{i}}
\title{$^{\mathcal{S}}$Linear operators in distribution spaces}
\date{}
\maketitle

\begin{abstract}
In this paper we define the $^{\mathcal{S}}$linear operators (Schwartz
linear operators) among spaces of tempered distributions. These operators
are the analogous of linear continuous operators among separable Hilbert
spaces, but in the case of spaces endowed with Schwartz bases having a
continuous index set. The Schwartz linear operatos enjoy properties very
similar to those enjoyed by linear operators in the finite dimensional case.
The $^{\mathcal{S}}$operators are one possible rigorous mathematical model
for the operators and observables used in Quantum Mechanics.
\end{abstract}

\bigskip 

\section{\textbf{Introduction}}

\bigskip

Let $X$ and $Y$ be two vector spaces on the field $\Bbb{K\;}$(the real field 
$\Bbb{R}\;$or the complex one $\Bbb{C}$). A function $f$ from $X$ into $Y$
is called linear if, for any two points $x,y$ of the space $X$ and for each
scalar $\lambda \in \Bbb{K}$, the equality 
\[
f(\lambda x+y)=\lambda f\left( x\right) +f(y), 
\]
holds true. Equivalently, a mapping $f$ from $X$ into $Y$ is linear if and
only if for every integer $k\in \Bbb{N}$, for any $k$-tuple $x=\left(
x_{i}\right) _{i=1}^{k}$ of points of the space $X$ and for any $k$-tuple of
scalars $\lambda =\left( \lambda _{i}\right) _{i=1}^{k}$ in $\Bbb{K}$,
setting 
\[
\sum \lambda x:=\sum_{i=1}^{k}\lambda _{i}x_{i} 
\]
and $f\left( x\right) :=\left( f\left( x_{i}\right) \right) _{i=1}^{k}$, we
have 
\[
f\left( \sum \lambda x\right) =\sum \lambda f\left( x\right) , 
\]
i.e., the image of the $\lambda $-linear combination of a family $x$ is the $%
\lambda $-linear combination of the image family $f(x)$ of the family $x$
under the function $f$; in indexed notation, we have 
\[
f\left( \sum_{i=1}^{k}\lambda _{i}x_{i}\right) =\sum_{i=1}^{k}\lambda
_{i}f\left( x_{i}\right) . 
\]
The aim of this chapter is to extend the last definition to the class of $^{%
\mathcal{S}}$families of tempered distributions indexed by the Euclidean
space $\Bbb{R}^{k}$, using, as coefficient systems, locally summable maps
from $\Bbb{R}^{k}$ to $\Bbb{K}$ and, more generally, Schwartz tempered
distributions from $\Bbb{R}^{k}$ into $\Bbb{K}$ (which, as we already have
seen, are so viewed as ``non-pointwise defined'' families in the field $\Bbb{%
K}$ indexed by the Euclidean space $\Bbb{R}^{k}$). If $v=\left( v_{i}\right)
_{i\in \Bbb{R}^{k}}$ is an $^{\mathcal{S}}$family in the distribution space $%
\mathcal{S}_{n}^{\prime }$, i.e. if for every test function $\phi \in 
\mathcal{S}_{n}$, the function 
\[
v\left( \phi \right) :\Bbb{R}^{k}\rightarrow \Bbb{K}:i\mapsto v_{i}\left(
\phi \right) , 
\]
belongs to the test function space $\mathcal{S}_{k}$, and if $\lambda \in 
\mathcal{S}_{k}^{\prime }$ is a tempered distribution defined on the index
set of the family $v$, we put 
\[
\int_{\Bbb{R}^{k}}\lambda v:=\lambda \circ \widehat{v}=\;^{t}\left( \widehat{%
v}\right) \left( \lambda \right) , 
\]
where $^{t}\left( \widehat{v}\right) $ is the (topological) transpose of the
continuous operator 
\[
\widehat{v}:\mathcal{S}_{n}\rightarrow \mathcal{S}_{k}:\phi \mapsto v\left(
\phi \right) . 
\]
The idea is very natural:

\begin{itemize}
\item  \emph{an operator }$L:\mathcal{S}_{n}^{\prime }\rightarrow \mathcal{S}%
_{m}^{\prime }$\emph{\ is said to be }$^{\mathcal{S}}$\emph{linear if, for
every integer }$k\in \Bbb{N}$\emph{, for every distribution coefficient }$%
\lambda $\emph{\ in }$\mathcal{S}_{k}^{\prime }$\emph{\ and for every family
of distributions }$v$\emph{\ in }$\mathcal{S}(\Bbb{R}^{k},\mathcal{S}%
_{n}^{\prime })$\emph{, the image of the }$^{\mathcal{S}}$\emph{family }$v$%
\emph{\ is an }$^{\mathcal{S}}$\emph{family too and the equality } 
\[
L\left( \int_{\Bbb{R}^{k}}\lambda v\right) =\int_{\Bbb{R}^{k}}\lambda
L\left( v\right) ,
\]
\emph{holds true.}
\end{itemize}

\bigskip

\section{$^{\mathcal{S}}$\textbf{Operators}}

\bigskip

First of all we have to transform a family of tempered distributions by
means of operators defined on spaces of tempered distributions, the
definition is pointwise and absolutely straightforward.

\bigskip

\textbf{Definition (image of a family of distributions).}\ \emph{Let }$W$%
\emph{\ be a subset of the distribution space }$\mathcal{S}_{n}^{\prime }$%
\emph{, let }$A:W\rightarrow \mathcal{S}_{m}^{\prime }$\emph{\ be an
operator (not necessarily linear) and let }$v=(v_{p})_{p\in \Bbb{R}^{k}}$%
\emph{\ be a family of tempered distributions belonging to the subset }$W$%
\emph{, i.e. a family with trace (trajectory) set }$\{v_{p}\}_{p\in \Bbb{R}%
^{k}}\emph{\ }$\emph{contained in the subset }$W$\emph{. The \textbf{image
of the family} }$v$ \emph{\textbf{by means of the operator} }$A$\emph{\ is,
by definition, the family }$A(v)$\emph{\ in }$\mathcal{S}_{m}^{\prime }$%
\emph{\ defined by} 
\[
A(v):=(A(v_{p}))_{p\in \Bbb{R}^{k}}, 
\]
\emph{i.e., the family }$A(v)$ \emph{such that, for all index }$p\in \Bbb{R}%
^{k}$\emph{, we have} $A(v)_{p}=A(v_{p})$\emph{.}

\bigskip

We can read the above definition saying that:

\begin{itemize}
\item  \emph{the image (under an operator) of a family of vectors is the
family of the images of vectors.}
\end{itemize}

\bigskip

\textbf{Definition (operator of class }$\mathcal{S}$\textbf{).}\emph{\ Let }$%
W$\emph{\ be a subset of the space }$\mathcal{S}_{n}^{\prime }$\emph{\ and
let }$L:W\rightarrow \mathcal{S}_{m}^{\prime }$\emph{\ be an operator (not
necessarily linear). The operator }$L$\emph{\ is said to be an }$^{\mathcal{S%
}}$\emph{\textbf{operator }or\textbf{\ operator of class }}$\mathcal{S}$%
\emph{\ if, for each natural }$k$\emph{\ and for each }$^{\mathcal{S}}$\emph{%
family }$v\in \mathcal{S}(\Bbb{R}^{k},\mathcal{S}_{n}^{\prime })$\emph{\
with trajectory contained in }$W\emph{,}$\emph{\ the image }$L(v)$\emph{\ of
the family }$v$ \emph{is an }$^{\mathcal{S}}$\emph{family too (that is, if
the image }$L(v)$\emph{\ belongs to the space }$\mathcal{S}(\Bbb{R}^{k},%
\mathcal{S}_{m}^{\prime })$\emph{).}

\bigskip

We can read the above definition as follows:

\begin{itemize}
\item  \emph{an operator }$L$\emph{\ is of class }$\mathcal{S}$\emph{\ if
the image by }$L$\emph{\ of any }$^{\mathcal{S}}$\emph{family is an }$^{%
\mathcal{S}}$\emph{family too.}
\end{itemize}

\bigskip

\section{$^{\mathcal{S}}$\textbf{Operators defined on }$\mathcal{S}
_{n}^{\prime }$}

\bigskip

The following property proves that the class of linear $^{\mathcal{S}}$%
operators defined on the entire space of tempered distribution contains the
class of weakly* continuous linear operators on that space.

\bigskip

\textbf{Theorem (the transpose of an operator).}\emph{\ The transpose of a
weakly continuous linear operator defined among two spaces of Schwartz test
functions is an }$^{\mathcal{S}}$\emph{operator. Consequently, every weakly*
continuous linear operator defined among two spaces of tempered
distributions is an }$^{\mathcal{S}}$\emph{operator. Moreover, the operator
associated with the image of a family }$v$\emph{\ by the transpose of a
weakly continuous operator }$A$\emph{\ is the composition }$\widehat{v}\circ
A$\emph{, that is, we have} 
\[
^{t}A\left( v\right) ^{\wedge }=\widehat{v}\circ A. 
\]

\emph{\bigskip }

\emph{Proof.} Let $A:\mathcal{S}_{n}\rightarrow \mathcal{S}_{m}$ be a
continuous linear operator with respect to the pair of weak topologies $%
\left( \sigma \left( \mathcal{S}_{n}\right) ,\sigma \left( \mathcal{S}%
_{m}\right) \right) $. Then, the operator $A$ is (topologically)
transposable (i.e., for every tempered distribution $a\in \mathcal{S}%
_{m}^{\prime }$, the functional $a\circ A$ lies in the space $\mathcal{S}%
_{n}^{\prime }$) and its (topological) transpose is (by definition) the
operator 
\[
^{t}A:\mathcal{S}_{m}^{\prime }\rightarrow \mathcal{S}_{n}^{\prime
}:a\mapsto a\circ A. 
\]
Let $v\in \mathcal{S}(\Bbb{R}^{k},\mathcal{S}_{n}^{\prime })$ be an $^{%
\mathcal{S}}$family of distributions, we have, by definition of image of a
family, 
\[
^{t}A\left( v\right) _{p}=\;^{t}A(v_{p}), 
\]
and hence we deduce 
\begin{eqnarray*}
^{t}A\left( v\right) \left( \phi \right) (p) &=&\;^{t}A\left( v\right)
_{p}\left( \phi \right) \\
&=&\;^{t}A(v_{p})\left( \phi \right) \\
&=&v_{p}\left( A\left( \phi \right) \right) \\
&=&v\left( A\left( \phi \right) \right) (p),
\end{eqnarray*}
so, taking into account that the family $v$ is an $^{\mathcal{S}}$family, we
deduce that the image 
\[
^{t}A\left( v\right) \left( \phi \right) =\widehat{v}\left( A\left( \phi
\right) \right) 
\]
belongs to the space $\mathcal{S}_{k}$. Concluding, the image family $%
^{t}A\left( v\right) $ is a family of class $\mathcal{S}$ belonging to the
space $\mathcal{S}(\Bbb{R}^{k},\mathcal{S}_{n}^{\prime })$, and thus the
transpose operator $^{t}A$, sending $^{\mathcal{S}}$families into $^{%
\mathcal{S}}$families, is an $^{\mathcal{S}}$operator. By the way, we proved
also that the operator associated with the image family $^{t}A(v)$ is the
composition $\widehat{v}\circ A$, that is $^{t}A\left( v\right) ^{\wedge }=%
\widehat{v}\circ A$. $\blacksquare $

\bigskip

\textbf{Application.} Let\textbf{\ }$L:\mathcal{S}_{n}^{\prime }\rightarrow 
\mathcal{S}_{n}^{\prime }$ be a differential operator with constant
coefficients and let $v$ be an $^{\mathcal{S}}$family in the space $\mathcal{%
S}_{n}^{\prime }$. Then $L(v)$ is an $^{\mathcal{S}}$family, in fact the
operator $L$ is the transpose of some differential operator on the space $%
\mathcal{S}_{n}$. For instance, the Dirac family $\left( \delta _{x}\right)
_{x\in \Bbb{R}^{n}}$ is obviously an $^{\mathcal{S}}$family, and so the
family of $i$-th derivatives $(\delta _{x}^{(i)})_{x\in \Bbb{R}^{n}}$ is an $%
^{\mathcal{S}}$family too, for every multi-index $i$.

\bigskip

\section{\textbf{Transposability of linear }$^{\mathcal{S}}$\textbf{%
operators on }$\mathcal{S}_{n}^{\prime }$}

\bigskip

The following property proves that linear $^{\mathcal{S}}$operators defined
among distribution spaces are weakly* continuous upon any subspace generated
(in the usual algebraic sense) by an $^{\mathcal{S}}$basis.

\bigskip

\textbf{Theorem (transposability of }$^{\mathcal{S}}$\textbf{operators). }%
\emph{A linear }$^{\mathcal{S}}$\emph{operator defined among two spaces of
tempered distributions is weakly* continuous on any algebraic linear hull of 
}$^{\mathcal{S}}$\emph{basis.}

\emph{\bigskip }

\emph{Proof. }If $L:\mathcal{S}_{n}^{\prime }\rightarrow \mathcal{S}%
_{m}^{\prime }$ is a linear $^{\mathcal{S}}$operator, and if $v$ is an $^{%
\mathcal{S}}$basis in $\mathcal{S}_{n}^{\prime }$ indexed by the $k$%
-dimensional real Euclidean space, then the image $L(v)$ of the $^{\mathcal{S%
}}$basis $v$ is a family of class $\mathcal{S}$ in the space $\mathcal{S}%
_{m}^{\prime }$ indexed by $\Bbb{R}^{k}$. So for every test function $h$ in $%
\mathcal{S}_{m}$, the image of the function $h$ by the family $L(v)$ is a
function of class $\mathcal{S}$ (belonging to the space $\mathcal{S}_{k}$),
namely the function 
\[
L(v)(h):\Bbb{R}^{k}\rightarrow \Bbb{K}:q\mapsto L(v_{q})(h). 
\]
We claim that the restriction $M$ of the operator $L$ to the pair $(E,F)$,
where $E$ is the linear hull of the family $v$ and $F$ is the linear hull of
the family $L(v)$, is weakly topologically transposable with respect to the
weak dual pair $(E,\mathcal{S}_{n})$ and $(F,\mathcal{S}_{m})$. Indeed, we
claim that its weak transpose $^{t}M$ is the operator 
\[
T:\mathcal{S}_{m}\rightarrow \mathcal{S}_{n}:h\mapsto v^{-}(L(v)(h)), 
\]
where $v^{-}$ is the inverse of the $^{\mathcal{S}}$basis $v$. We have to
prove, by the classic definition of weak transpose, that 
\[
\left\langle u,T(h)\right\rangle _{n}=\left\langle L(u),h\right\rangle _{m}, 
\]
for every tempered distribution $u$ in the space $E$ and for every test
function $h$ in the space $\mathcal{S}_{m}$. But, since $u$ is a finite
linear combination of the family $v$, we can prove the above duality
condition only for the elements of the family $v$, and we indeed have, for
every $q$ in $\Bbb{R}^{k}$, 
\begin{eqnarray*}
\left\langle v_{q},T(h)\right\rangle _{n} &=&\left\langle
v_{q},v^{-}(L(v)(h))\right\rangle _{n}= \\
&=&v(v^{-}(L(v)(h)))(q)= \\
&=&L(v)(h)(q)= \\
&=&L(v)_{q}(h)= \\
&=&L(v_{q})(h)= \\
&=&\left\langle L(v_{q}),h\right\rangle _{m},
\end{eqnarray*}
as we desire. Since every weak topologically transposable operator is weakly
continuous, we conclude that every linear $^{\mathcal{S}}$operator is weakly
continuous on the linear span of an $^{\mathcal{S}}$basis. $\blacksquare $

\bigskip

\section{$^{\mathcal{S}}$\textbf{Linear operators on }$\mathcal{S}
_{n}^{\prime }$}

\bigskip

In this section we shall introduce the main concept of the chapter.

\bigskip

\textbf{Definition (}$^{\mathcal{S}}$\textbf{linear operators on the entire }%
$\mathcal{S}_{n}^{\prime }$\textbf{).\ }\emph{Let} $L:\mathcal{S}%
_{n}^{\prime }\rightarrow \mathcal{S}_{m}^{\prime }$ \emph{be an }$^{%
\mathcal{S}}$\emph{operator (not necessarily linear). The operator }$L$\emph{%
\ is called }$^{\mathcal{S}}$\emph{\textbf{linear operator} if, for each
positive integer }$k$\emph{, for each }$^{\mathcal{S}}$\emph{family }$v\in 
\mathcal{S}(\Bbb{R}^{k},\mathcal{S}_{n}^{\prime })$\emph{\ and for every
tempered distribution }$a$\emph{\ in the space }$\mathcal{S}_{k}^{\prime }$%
\emph{, the equality} 
\[
L\left( \int_{\Bbb{R}^{k}}av\right) =\int_{\Bbb{R}^{k}}aL(v) 
\]
\emph{holds true.}

\bigskip

Utterly, an $^{\mathcal{S}}$linear operator must be linear, as we prove
below.

\bigskip

\textbf{Property (linearity of the }$^{\mathcal{S}}$\textbf{linear
operators).} \emph{An }$^{\mathcal{S}}$\emph{linear operator is linear.}

\emph{\bigskip }

\emph{Proof.} Indeed, in the conditions of the above definition, for each
couple of scalars $b,c$ and any couple $u,w$ of tempered distributions in
the space $\mathcal{S}_{n}^{\prime }$, if $\delta $ is the Dirac basis of
the space $\mathcal{S}_{n}^{\prime }$, we have 
\begin{eqnarray*}
L(au+bw) &=&L\left( \int_{\Bbb{R}^{k}}(au+bw)\delta \right) = \\
&=&\int_{\Bbb{R}^{k}}(au+bw)L(\delta )= \\
&=&a\int_{\Bbb{R}^{k}}uL(\delta )+b\int_{\Bbb{R}^{k}}wL(\delta )= \\
&=&aL\left( \int_{\Bbb{R}^{k}}u\delta \right) +bL\left( \int_{\Bbb{R}%
^{k}}w\delta \right) = \\
&=&aL(u)+bL(w),
\end{eqnarray*}
as we desired. $\blacksquare $

\bigskip

But we will see more than this preliminary remark about $^{\mathcal{S}}$%
linear operators.

\bigskip

\textbf{Remark.} The above definition and property can be immediately be
generalized to the case of operators defined on the $^{\mathcal{S}}$linear
hull $E$ of an $^{\mathcal{S}}$basis.

\bigskip

\section{\textbf{Examples of }$^{\mathcal{S}}$\textbf{linear operators}}

\bigskip

In this section we propose two important examples of $^{\mathcal{S}}$linear
operators. We note that the first is a particular case of the second one,
and indeed we shall see that every $^{\mathcal{S}}$linear operator defined
on the entire $\mathcal{S}_{n}^{\prime }$ is of the type presented in the
second example.

\bigskip

\subsection{\textbf{The superposition operator of an }$^{\mathcal{S}}$%
\textbf{family}}

\bigskip

Recall that if $v\in s(\Bbb{R}^{k},\mathcal{S}_{m}^{\prime })$ is any family
of tempered distributions\ in $\mathcal{S}_{m}^{\prime }$ indexed by an
Euclidean space $\Bbb{R}^{k}$ and if $w\in \mathcal{S}(\Bbb{R}^{m},\mathcal{S%
}_{n}^{\prime })$ is any $^{\mathcal{S}}$family of tempered distributions in 
$\mathcal{S}_{n}^{\prime }$, the family in $\mathcal{S}_{n}^{\prime }$
indexed by $\Bbb{R}^{k}$ and defined by 
\[
\int_{\Bbb{R}^{m}}vw:=\left( \int_{\Bbb{R}^{m}}v_{p}w\right) _{p\in \Bbb{R}%
^{k}}, 
\]
is called the superposition of the $^{\mathcal{S}}$family $w$\ with respect
to the family $v$.

\bigskip

We have already proved that, if the family $v$ belongs to the space $%
\mathcal{S}(\Bbb{R}^{k},\mathcal{S}_{m}^{\prime })$ then the superposition $%
\int_{\Bbb{R}^{m}}vw$ belongs to the space $\mathcal{S}(\Bbb{R}^{k},\mathcal{%
S}_{n}^{\prime })$ and the operator associated with this superposition is
the composition of the operators associated with the two families $v$ and $w$
, precisely we have 
\[
\left( \int_{\Bbb{R}^{m}}vw\right) ^{\wedge }=\widehat{v}\circ \widehat{w}. 
\]
In this case, sometimes, it is also convenient to denote the superposition 
\[
\int_{\Bbb{R}^{m}}vw 
\]
by the product notation $v.w$ and to call it also the $^{\mathcal{S}}$%
product of the family $v$ by the family $w$.

\bigskip

\textbf{Proposition.\ }\emph{Let }$w\in \mathcal{S}(\Bbb{R}^{m},\mathcal{S}%
_{n}^{\prime })$\emph{\ be an }$^{\mathcal{S}}$\emph{family of distributions
and let }$L:\mathcal{S}_{m}^{\prime }\rightarrow \mathcal{S}_{n}^{\prime }$%
\emph{\ be the superposition operator of the family }$w$\emph{, defined by} 
\[
L(a)=\int_{\Bbb{R}^{m}}aw, 
\]
\emph{for\ all tempered distribution }$a\in \mathcal{S}_{m}^{\prime }.$\emph{%
\ Then, the operator }$L$\emph{\ is an }$^{\mathcal{S}}$\emph{linear
operator.}

\emph{\bigskip }

\emph{Proof.} The operator $L$ is an $^{\mathcal{S}}$operator, indeed we
know that $L$ is the transpose of the continuous linear operator associated
with $v$ and then it is weakly* continuous. But we decide to see this fact
directly too. If $v\in \mathcal{S}(\Bbb{R}^{k},\mathcal{S}_{m}^{\prime })$
is an $^{\mathcal{S}}$family then its image by the operator $L$ is $L\left(
v\right) =v.w$ and the product of two $^{\mathcal{S}}$families is an $^{%
\mathcal{S}}$family. Let $a\in \mathcal{S}_{k}^{\prime }$ be a tempered
distribution and let $v\in \mathcal{S}(\Bbb{R}^{k},\mathcal{S}_{m}^{\prime
}) $ be an $^{\mathcal{S}}$family, we have 
\begin{eqnarray*}
L\left( \int_{\Bbb{R}^{k}}av\right) &=&\int_{\Bbb{R}^{m}}\left( \int_{\Bbb{R}%
^{k}}av\right) w= \\
&=&\int_{\Bbb{R}^{k}}a\left( \int_{\Bbb{R}^{m}}vw\right) = \\
&=&\int_{\Bbb{R}^{k}}aL(v),
\end{eqnarray*}
applying the already known property of ``$^{\mathcal{S}}$linearity'' of
superpositions. Note in fact that, for each index $p\in \Bbb{R}^{k}$, we
have 
\begin{eqnarray*}
L(v)_{p} &=&L(v_{p})= \\
&=&\int_{\Bbb{R}^{m}}v_{p}w= \\
&=&\left( \int_{\Bbb{R}^{m}}vw\right) _{p},
\end{eqnarray*}
and the proof is completed. $\blacksquare $

\bigskip

\subsection{\textbf{Transpose operators}}

\bigskip

\textbf{Lemma} \textbf{(the image under a transpose operator).}\emph{\ Let }$%
B\in \mathcal{L}\left( \mathcal{S}_{n},\mathcal{S}_{m}\right) $\emph{\ be a
linear continuous operator\ and let }$v\in \mathcal{S}(\Bbb{R}^{k},\mathcal{S%
}_{m}^{\prime })$\emph{\ be an }$^{\mathcal{S}}$\emph{family.\ Then, the
image of the family }$v$\emph{\ by the transpose operator }$^{t}B$\emph{\ is
the product of the family }$v$\emph{\ by the family generated by the
operator }$B$\emph{, in symbol we have} 
\[
^{t}B(v)=\int_{\Bbb{R}^{k}}vB^{\vee }, 
\]
\emph{so in particular, the transpose operator }$^{t}B$\emph{\ is an }$^{%
\mathcal{S}}$\emph{operator.}

\bigskip

\emph{Proof.} For each index $p\in \Bbb{R}^{k}$, we have 
\begin{eqnarray*}
\left( \int_{\Bbb{R}^{m}}vB^{\vee }\right) _{p} &=&\int_{\Bbb{R}
^{m}}v_{p}B^{\vee }= \\
&=&v_{p}\circ (B^{\vee })^{\wedge }= \\
&=&v_{p}\circ B= \\
&=&\;^{t}B(v_{p})= \\
&=&\;^{t}B(v)(p),
\end{eqnarray*}
and hence 
\[
\int_{\Bbb{R}^{m}}vB^{\vee }=\;^{t}B(v), 
\]
as we desired. $\blacksquare $

\bigskip

\textbf{Theorem} \textbf{(}$^{\mathcal{S}}$\textbf{linearity of a transpose
operator).} \emph{Let }$B\in \mathcal{L}(\mathcal{S}_{n},\mathcal{S}_{m})$%
\emph{\ be a linear and continuous operator and let }$v\in \mathcal{S}(\Bbb{R%
}^{k},\mathcal{S}_{m}^{\prime })$\emph{\ be an }$^{\mathcal{S}}$\emph{
family.\ Then, for each tempered coefficient system }$a\in \mathcal{S}%
_{k}^{\prime }$\emph{, we have} 
\[
^{t}B\left( \int_{\Bbb{R}^{k}}av\right) =\int_{\Bbb{R}^{k}}a\;^{t}B(v). 
\]

\bigskip

\emph{Proof.} We have 
\begin{eqnarray*}
^{t}B\left( \int_{\Bbb{R}^{k}}av\right) &=&\left( \int_{\Bbb{R}
^{k}}av\right) \circ B= \\
&=&(a\circ \widehat{v})\circ B= \\
&=&a\circ (\widehat{v}\circ B)= \\
&=&\int_{\Bbb{R}^{k}}a(\widehat{v}\circ B)^{\vee }= \\
&=&\int_{\Bbb{R}^{k}}a(\int_{\Bbb{R}^{m}}vB^{\vee })= \\
&=&\int_{\Bbb{R}^{k}}a\;^{t}B(v),
\end{eqnarray*}
as we desired. $\blacksquare $

\bigskip

\textbf{Application (derivatives of a distribution). }As a simple
application, we prove the formula 
\[
u^{\prime }=\int_{\Bbb{R}}u\delta ^{\prime }, 
\]
where $\delta ^{\prime }$ is the $^{\mathcal{S}}$family in $\mathcal{S}%
_{1}^{\prime }$ defined by $\delta ^{\prime }=(\delta _{p}^{\prime })_{p\in 
\Bbb{R}}$. Recall that the differential operators on the space of tempered
distributions are transpose of linear continuous operators and then they are 
$^{\mathcal{S}}$linear operators. Let $\delta $ be the Dirac family of the
space $\mathcal{S}_{1}^{\prime }$, then for each tempered distribution $u\in 
\mathcal{S}_{1}^{\prime }$, we have 
\[
u=\int_{\Bbb{R}}u\delta , 
\]
and consequently 
\begin{eqnarray*}
u^{\prime } &=&\partial \left( \int_{\Bbb{R}}u\delta \right) = \\
&=&\int_{\Bbb{R}}u\partial (\delta )= \\
&=&\int_{\Bbb{R}}u\delta ^{\prime }.
\end{eqnarray*}
More generally, in the space $\mathcal{S}_{n}^{\prime }$, we have (by the
same proof) 
\[
L(u)=\int_{\Bbb{R}}uL(\delta ), 
\]
for every differential operator $L$, and every tempered distribution $u$.

\bigskip

\section{\textbf{Characterization of }$^{\mathcal{S}}$\textbf{linear
operators}}

\bigskip

Now, we can show the true nature of the $^{\mathcal{S}}$linear operators
defined on $\mathcal{S}_{n}^{\prime }$.

\bigskip

\textbf{Theorem (characterization of }$^{\mathcal{S}}$\textbf{linearity).} 
\emph{Let }$L:\mathcal{S}_{n}^{\prime }\rightarrow \mathcal{S}_{m}^{\prime }$%
\emph{\ be an operator. Then, }$L$\emph{\ is }$^{\mathcal{S}}$\emph{linear
if and only if there exists a linear and continuous operator }$B\in \mathcal{%
L}\left( \mathcal{S}_{m},\mathcal{S}_{n}\right) $\emph{\ such that }$%
L=\;^{t}\left( B\right) $\emph{.}

\bigskip

\emph{Proof.} \emph{Sufficiency.} It follows from the above theorem. \emph{%
Necessity. }Let $\delta $ be the Dirac family in the space $\mathcal{S}%
_{n}^{\prime }$, we have 
\begin{eqnarray*}
L\left( u\right) &=&L\left( \int_{\Bbb{R}^{n}}u\delta \right) = \\
&=&\int_{\Bbb{R}^{n}}uL\left( \delta \right) = \\
&=&\;^{t}\left( L\left( \delta \right) ^{\wedge }\right) \left( u\right) ,
\end{eqnarray*}
so the operator $L$ is the transpose of the operator generated by the
Schwartz family $L(\delta )$, that is 
\[
L=\;^{t}\left( L\left( \delta \right) ^{\wedge }\right) , 
\]
since the operators generated by Schwartz families are continuous, we
conclude the proof. $\blacksquare $

\bigskip

Before to give the last complete characterization of $^{\mathcal{S}}$linear
operators, we recall the following classical definition from Linear
Functional Analysis.

\bigskip

\textbf{Definition (of transposable operator).} \emph{A linear operator }$L:%
\mathcal{S}_{n}^{\prime }\rightarrow \mathcal{S}_{m}^{\prime }$\emph{\ is
said to be \textbf{transposable with respect to the canonical pairings} }$(%
\mathcal{S}_{n},\mathcal{S}_{n}^{\prime })$\emph{\ \textbf{and} }$(\mathcal{S%
}_{m},\mathcal{S}_{m}^{\prime })$\emph{\ if and only if there exists a
linear continuous operator }$B\in \mathcal{L}(\mathcal{S}_{m},\mathcal{S}
_{n})$\emph{\ such that }$L=\;^{t}\left( B\right) $\emph{.}

\bigskip

Recalling that the operator $L$ is weakly continuous if and only if it is
strongly continuous if and only if it is transposable, we derive the
following definitive characterization.

\bigskip

\textbf{Theorem (characterization of }$^{\mathcal{S}}$\textbf{linearity).} 
\emph{Let }$L:\mathcal{S}_{n}^{\prime }\rightarrow \mathcal{S}_{m}^{\prime }$%
\emph{\ be an operator. Then, the following assertions are equivalent}

\begin{itemize}
\item[\emph{1)}]  \emph{\ the operator }$L$\emph{\ is }$^{\mathcal{S}}$\emph{%
linear;}

\item[\emph{2)}]  \emph{\ there exists an operator }$B\in \mathcal{L}\left( 
\mathcal{S}_{m},\mathcal{S}_{n}\right) $\emph{\ such that }$L=\;^{t}\left(
B\right) $\emph{;}

\item[\emph{3)}]  \emph{\ the operator }$L$\emph{\ is linear and weakly
continuous;}

\item[\emph{4)}]  \emph{\ the operator }$L$\emph{\ is linear and strongly
continuous;}

\item[\emph{5)}]  \emph{\ the operator }$L$\emph{\ is linear and
topologically transposable.}
\end{itemize}

$\bigskip $

$\bigskip $

\bigskip

\textbf{David Carf\`{i}}

\emph{Faculty of Economics}

\emph{University of Messina}

\emph{davidcarfi71@yahoo.it}

\end{document}